\documentclass[10pt,reqno]{amsart}
\usepackage[T1]{fontenc}
\usepackage{lmodern}
\usepackage[utf8]{inputenc}
\usepackage[english]{babel}
\usepackage{amsmath,amssymb,mathtools}
\usepackage[shortlabels]{enumitem}
\usepackage[autostyle]{csquotes}
\usepackage[colorlinks=true,linkcolor=blue,citecolor=blue,urlcolor=blue]{hyperref}
\usepackage{orcidlink}
\usepackage[backend=bibtex,sorting=nyt,style=alphabetic,maxbibnames=9]{biblatex}
\allowdisplaybreaks
\title[IMC and NSOP$_1$]
{Idempotent generically stable measures in NSOP$_1$ groups}
\author[D. M. Hoffmann]{Daniel Max Hoffmann$^{\dagger}$\
\orcidlink{0000-0002-4514-269X}}
\thanks{$^{\dagger}$SDG. The author is supported by the National Science Centre
(Narodowe Centrum Nauki, Poland) grant no.~2021/43/B/ST1/00405,
and by the UW IDUB PSP no.~501-D110-20-0001381.}
\address{$^{\dagger}$Instytut Matematyki\\Uniwersytet Warszawski\\Warszawa\\Poland}
\email{daniel.max.hoffmann@gmail.com}
\urladdr{https://sites.google.com/site/danielmaxhoffmann/}
\date{November 23, 2026}

\newtheorem{theorem}{Theorem}[section]
\newtheorem{proposition}[theorem]{Proposition}
\newtheorem{lemma}[theorem]{Lemma}
\newtheorem{cor}[theorem]{Corollary}
\newtheorem{fact}[theorem]{Fact}
\theoremstyle{definition}
\newtheorem{definition}[theorem]{Definition}

\newtheorem*{theorem*}{Theorem}
\newcommand{\FC}{\mathfrak C}
\newcommand{\CL}{\mathcal L}
\newcommand{\NSOP}{\mathrm{NSOP}_1}
\newcommand{\SOP}{\mathrm{SOP}_1}
\newcommand{\checkmu}{\check\mu}
\DeclareMathOperator{\supp}{supp}
\DeclareMathOperator{\stab}{Stab}
\DeclareMathOperator{\tp}{tp}

\DeclareMathOperator{\E}{\mathbb E}

\newtheorem{clm}{Claim}
\newtheorem*{clm*}{Claim}

\AtEndEnvironment{proof}{\setcounter{clm}{0}}
\newenvironment{clmproof}[1][\proofname]{\proof[#1]}{\endproof}

\def\Ind#1#2{#1\setbox0=\hbox{$#1x$}\kern\wd0\hbox to 0pt{\hss$#1\mid$\hss}
\lower.9\ht0\hbox to 0pt{\hss$#1\smile$\hss}\kern\wd0}

\def\ind{\mathop{\mathpalette\Ind{}}}

\def\notind#1#2{#1\setbox0=\hbox{$#1x$}\kern\wd0
\hbox to 0pt{\mathchardef\nn=12854\hss$#1\nn$\kern1.4\wd0\hss}
\hbox to 0pt{\hss$#1\mid$\hss}\lower.9\ht0 \hbox to 0pt{\hss$#1\smile$\hss}\kern\wd0}

\begin{document}
\begin{abstract}
We prove the Idempotent Measure Conjecture for type-definable groups in theories
whose Keisler randomization is NSOP$_1$. More precisely, every global
idempotent frequency interpretation measure on such a group is the unique
left-invariant and the unique right-invariant measure concentrating on its
type-definable stabilizer. In particular, its support is a left-simple
convolution semigroup and has a unique minimal left ideal.

The announced theorem of Ben Yaacov, Chernikov and
Ramsey that randomization preserves NSOP$_1$ (proof in the appendix)
yields the corresponding
conclusion for all type-definable groups in NSOP$_1$ theories.
\end{abstract}
\maketitle

\section{Introduction}
The classification of idempotent probability measures relates the
algebra of convolution to the structure of subgroups. In the classical
setting, Pym's theorem states that a regular Borel probability measure on
a locally compact group is idempotent if and only if it is the normalized
Haar measure of a compact subgroup \cite[Theorem~A.4.1]{Pym62}.
The Idempotent Measure Conjecture asks for an analogous characterization
for Keisler measures on type-definable groups, with type-definable subgroups
in place of compact subgroups (cf. \cite{Artem_Kyle}, \cite{Artem_Kyle2}, \cite{CGK}).

Let $G$ be a group type-definable in a monster model $\FC$ of a complete theory
$T$. Definable convolution is the pushforward of the Morley product under
group multiplication. We consider frequency interpretation measures (\emph{fim measures}) introduced in \cite{HruPiSi13}. In NIP theories these
are precisely the generically stable measures; outside NIP we use the fim
condition throughout. For a global fim measure $\mu$, the subgroup
\[
 H_\mu=\stab_{\mathrm r}(\mu)
      =\{g\in G(\FC):\mu*\delta_g=\mu\}
\]
is type-definable (cf. \cite[Fact~3.35]{CGK}). An idempotent fim measure is called
\emph{generically transitive} if it concentrates on $H_\mu$.
The conjecture asks whether every global fim idempotent is generically
transitive \cite[Problem~3.41]{CGK}.

Chernikov and Gannon proved the conjecture in stable theories
\cite[Theorem~5.8]{Artem_Kyle}. Chernikov, Gannon and Krupi\'nski proved
it for abelian groups in arbitrary theories \cite[Theorem~3.45]{CGK},
and the result was extended to virtually solvable definable groups in
\cite{HoffmannSOLVABLE}. Quite recently, Chernikov resolved the conjecture for
all type-definable groups in NIP theories \cite{ChernikovNIP2026}.
The present paper treats a different hypothesis on the ambient theory.
Write $T^R$ for its Keisler randomization, viewed as a continuous theory.

\begin{theorem}\label{thm:main}
Suppose that $T^R$ is $\NSOP$, let $G$ be a type-definable group in $T$, and
let $\mu\in\mathfrak M_G(\FC)$ be a global fim measure. Then the following
conditions are equivalent:
\begin{enumerate}[(i)]
\item $\mu*\mu=\mu$;
\item $\mu$ concentrates on $H=\stab_{\mathrm r}(\mu)$.
\end{enumerate}
When these conditions hold, $H=\stab_{\ell}(\mu)$ and $\checkmu=\mu$,
where $\checkmu$ is the pushforward under inversion. Among global measures
concentrating on $H$, the measure $\mu$ is the unique left-$H$-invariant
measure and the unique right-$H$-invariant measure. In particular, $H$ is
a fim group and has finitely satisfiable generics.
\end{theorem}

Here, a type-definable group $H$ is called a \emph{fim group}
if it admits a global right-$H$-invariant fim measure
concentrating on $H$. It has \emph{finitely satisfiable generics}
(\emph{fsg}) if there are a small model $M$ over which $H$ is
type-definable and a global type $p\in S_H(\FC)$ such that
every left translate $g\cdot p$, for $g\in H(\FC)$, is finitely
satisfiable in $H(M)$; see
\cite[Definitions~3.27 and~3.29]{CGK}.

Ben Yaacov, Chernikov and Ramsey have announced that
randomization preserves $\NSOP$; see the seminar announcement
\cite{BenYaacovNSOP2019} and the slides \cite{ChernikovFields2021}.
We could not find the proof, so we provide one in the appendix.
Using this result gives the following consequence.

\begin{cor}\label{cor:ground-main}
Every global fim idempotent measure on a type-definable group
in an $\NSOP$ theory is generically transitive, and satisfies all the
conclusions of Theorem~\ref{thm:main}.
\end{cor}

Our main technical result is the following obstruction theorem
($T$ is arbitrary):

\begin{proposition}\label{thm:obstruction-intro}
If a type-definable group in $T$ carries a global fim idempotent measure $\nu$
such that $\nu\ne\check\nu$, then $T^R$ has $\SOP$.
\end{proposition}

To deduce generic transitivity, we use the identity
$\mu*\delta_b*\mu=\mu$ for $b$ realizing a type in $\supp(\mu|M)$, where $M$ is a small
model over which $\mu$ is fim \cite[Lemma~3.1]{HoffmannSOLVABLE} (see also Lemma 3.4 in \cite{ChernikovNIP2026}).
Both $\mu$ and $\mu*\delta_b$ are then fim idempotents. Their invariance
under inversion forces $\mu*\delta_b=\mu$. This last algebraic step is
also used in the proof of the NIP theorem in \cite[Section~5]{ChernikovNIP2026}.
Our proof of inversion invariance uses a tree configuration rather than
the bounded alternation argument of that paper.

Section~\ref{sec:preliminaries} collects the measure-theoretic facts and
the reduction to inversion invariance. Section~\ref{sec:comb} gives the
continuous tree criterion. Section~\ref{sec:asymmetry} proves the
obstruction result, the main theorem and records some of its consequences.
The randomization preservation theorem of Ben Yaacov, Chernikov and Ramsey is discussed in the Appendix, Section \ref{app:preservation}.

\ \\ 
\textbf{Use of AI.} The author used ChatGPT (GPT-6 Astra)
to assist in developing proofs and drafting the text.
Except for the appendix, 
the results were obtained by mid-November
and then polished.
Most of the work on the appendix was undertaken by the AI.
The author subsequently refined and streamlined the arguments and assumes full responsibility for the final content.

\section{Preliminaries}
\label{sec:preliminaries}

We work in a complete first-order theory $T$ with a monster model $\FC$.
All parameter sets, elementary submodels and tuples of variables
are small unless stated otherwise. 
Groups are type-definable in small tuples, with
type-definable graphs of multiplication and inversion. 
Writing $G\subseteq\FC^x$, by compactness we may fix maps
\[
 m:\FC^x\times\FC^x\longrightarrow\FC^x,
 \qquad \iota:\FC^x\longrightarrow\FC^x
\]
extending multiplication and inversion on $G$, respectively,
such that each output coordinate is a definable function of
finitely many input coordinates. We use these extensions when
evaluating group words on ambient tuples. 

All parameter sets used as bases contain the parameters defining
the groups and the chosen extensions. We write $\CL(M)$ for formulas
over $M$.

\subsection{Fim measures and convolution}
For a tuple $x$ and a parameter set $A$, let $\mathfrak M_x(A)$ denote
the space of Keisler measures over $A$ in variables $x$. We identify a
measure with its regular Borel probability measure on $S_x(A)$.
Its support consists of the types all of whose formulas have positive
measure. A global measure is a measure over $\FC$.
For a type-definable set $P$, concentration on $P$ means that every
definable set containing $P$ has measure one.
For $f\in C(S_x(M))$, we write
$\mu(f)=\int f\,d(\mu|M)$.

A global measure $\mu_x$ is \emph{definable over $M$} if, for every
formula $\varphi(x;y)$, the function
\[
 S_y(M)\longrightarrow[0,1],\qquad
 \tp(b/M)\longmapsto\mu(\varphi(x;b))
\]
is well-defined and continuous. If $\mu$ is definable and $\nu_y$
is any global measure, their Morley product is given by
\[
 (\mu_x\otimes\nu_y)(\theta(x,y;c))
 =\int_{S_y(N)}\mu(\theta(x;b;c))\,d(\nu|N)(\tp(b/N)),
\]
where $N$ contains $c$ and a base of definability of $\mu$.

\begin{definition}\label{def:fim}
A global measure $\mu_x$ is \emph{fim over $M$} if it is definable over
$M$ and, for every formula $\varphi(x;y)$, there are formulas
$\theta_n(x_1,\ldots,x_n)$ over $M$ such that
$\mu^{(n)}(\theta_n)\to1$ and, for every $\varepsilon>0$ and all
sufficiently large $n$,
\[
 \sup_b\left|\frac1n\sum_{i=1}^n\mathbf1_{\varphi(a_i;b)}
                  -\mu(\varphi(x;b))\right|<\varepsilon
 \quad\text{whenever }\models\theta_n(a_1,\ldots,a_n).
\]
Here $\mu^{(n)}$ is the $n$-fold Morley product of $\mu$.
\end{definition}

We use the following standard properties
\cite[Section~5]{InTheWild}, \cite[Sections~3.3 and~3.5]{CGK}.
Definable pushforwards preserve the fim property, after enlarging the
base to include their parameters. Fim measures commute with one another
under the Morley product. Their finite Morley powers are therefore
invariant under permutations of the coordinates. Morley products of
definable measures are associative and commute with coordinatewise
definable pushforwards.

For definable measures on $G$, put
\[
 \mu*\nu=m_*(\mu\otimes\nu),\qquad m(x,y)=xy,
 \qquad \check\mu=(x\mapsto x^{-1})_*\mu.
\]
Convolution is associative. Commutation of Morley products implies, for
fim $\mu,\nu$, the identity
\begin{equation}\label{eq:inverse-convolution}
 \check{(\mu*\nu)}=\check\nu*\check\mu.
\end{equation}
Indeed, inversion sends $xy$ to $y^{-1}x^{-1}$, and the Morley factors
may be interchanged before applying this word map.
In particular, the pushforward of a fim idempotent under inversion is again
a fim idempotent.

The following block-product identity is also recorded in
\cite[Remark~2.5(iii)]{ChernikovNIP2026}.

\begin{lemma}\label{lem:blocks}
Let $\mu$ be a fim idempotent. Suppose that $W_1,\ldots,W_s$ are
ordered products of distinct variables from $x_1,\ldots,x_n$, each
nonempty, and that no variable occurs in two of these products. Then
\[
 ((x_1,\ldots,x_n)\mapsto(W_1,\ldots,W_s))_*\mu^{(n)}
 =\mu^{(s)}.
\]
\end{lemma}
\begin{proof}
Clear by idempotence and symmetry of Morley product of fim measures (\cite[Theorem 5.16(a)]{InTheWild}).
\end{proof}

\subsection{Generic transitivity}
For a global definable measure $\mu$ on $G$, its left and right
stabilizers are
\[
 \stab_{\ell}(\mu)=\{g\in G(\FC):\delta_g*\mu=\mu\},\qquad
 \stab_{\mathrm r}(\mu)=\{g\in G(\FC):\mu*\delta_g=\mu\}.
\]
They are type-definable over a base of definability of $\mu$
\cite[Fact~3.35]{CGK}. A global fim idempotent is
\emph{generically transitive} if it concentrates on its right stabilizer.
Equivalently,
\begin{equation}\label{eq:shear}
 ((x,y)\mapsto(xy,y))_*\mu^{(2)}=\mu^{(2)}
 \qquad\text{\cite[Proposition~3.37]{CGK}.}
\end{equation}

\begin{fact}\label{fact:classification}
Let $\mu$ be a global fim measure concentrating on
$H=\stab_{\mathrm r}(\mu)$. Then $\mu$ is idempotent,
$H=\stab_{\ell}(\mu)$, and $\mu$ is the unique left-$H$-invariant
and the unique right-$H$-invariant global measure concentrating on $H$.
Moreover, $\checkmu=\mu$, and $H$ is a fim group, hence an fsg group.
\end{fact}

 \begin{proof}
Follows by \cite[Remark~3.30, Proposition~3.33 and Remark~3.39]{CGK}.
\end{proof}

\begin{fact}\label{fact:sandwich}
Let $\mu$ be a global fim idempotent over $M$, on an $M$-type-definable group
$G$. If $\tp(b/M)\in\supp(\mu|M)$, then
\begin{equation}\label{eq:sandwich}
 \mu*\delta_b*\mu=\mu.
\end{equation}
\end{fact}
\begin{proof}
The proof of \cite[Lemma~3.1]{HoffmannSOLVABLE} applies equally to
type-definable groups.
\end{proof}

\begin{proposition}\label{prop:symmetry-reduction}
Fix a type-definable group $G$ in an arbitrary theory. The following are
equivalent:
\begin{enumerate}[(i)]
\item every global fim idempotent on $G$ is generically transitive;
\item every global fim idempotent on $G$ is invariant under inversion.
\end{enumerate}
\end{proposition}
\begin{proof}
The forward implication is Fact~\ref{fact:classification}.
Assume (ii), and let $\mu$ be a fim idempotent over $M$. Fix
$b\in G(\FC)$ with $\tp(b/M)\in\supp(\mu|M)$ and put
$\rho=\mu*\delta_b$. This is fim over any small model containing
$Mb$. By \eqref{eq:sandwich},
\begin{equation}\label{eq:rectangular}
 \mu*\rho=\rho,\qquad \rho*\mu=\mu,\qquad \rho*\rho=\rho.
\end{equation}
By (ii), both $\mu$ and $\rho$ are invariant under inversion.
Inverting the first equality in \eqref{eq:rectangular}, and using
\eqref{eq:inverse-convolution}, gives $\rho*\mu=\rho$.
The second equality yields $\rho=\mu$. Thus every realization of
$\supp(\mu|M)$ belongs to $\stab_{\mathrm r}(\mu)$.
Since this subgroup is $M$-type-definable, $\mu$ concentrates on it.
\end{proof}

\subsection{Keisler randomization}
We use Keisler randomization in the form of
\cite{BenYaacov2009}. Let $\mathcal R$ be a sufficiently saturated
model of $T^R$, containing a randomized elementary extension of $M$.
Write $M^c$ for the copy of $M$, consisting of constant
random variables. If $\varphi$ is a first-order formula, its expectation
in the randomization is denoted by $\E[\varphi]$.

A Keisler measure $\lambda$ over $M$ in a finite tuple $\bar x$
determines a unique randomized type over $M^c$ by
\[
 \E[\varphi(\bar X;m)]=\lambda(\varphi(\bar x;m))
 \qquad(\varphi(\bar x;m)\in\CL(M)).
\]
Existence follows by realizing the probabilities assigned by $\lambda$ to the
Boolean combinations of finitely many formulas on an atomless
probability space, and then applying compactness. Uniqueness follows from
quantifier elimination for randomization
\cite[Theorem~2.9 and Corollary~2.10]{BenYaacov2009}.
The same statement for any small tuple follows by compactness from its
finite subtuples. We call $\lambda$ the \emph{law over $M$} of $\bar X$.
Definable maps act on laws by pushforward.

Let $G=\pi(\FC)\subseteq\FC^x$, with $\pi$ over $M$. Write
$G^R$ for the set of realizations of the partial type
\[
 \pi^R(X)=\{\E[\theta(X)]=1:\theta\in\pi\}.
\]
The ambient definable extensions of multiplication and inversion
preserve $G^R$ and satisfy the group identities there.
We extend every finite-valued
relatively $M$-definable function $h$ on $G$ to the ambient domain,
preserving its range, and use the same notation for the extension.
Thus $\E[h(X)]$ is a continuous
formula in randomized variables, obtained as a finite
linear combination of expressions $\E[\varphi(X)]$, where
$\varphi$ is a first-order formula.

\section{A tree criterion for continuous \texorpdfstring{$\mathrm{SOP}_1$}{SOP1}}
\label{sec:comb}

We use the following convention for $\mathrm{SOP}_1$ in continuous logic.
A continuous formula $\varphi(x;y)$ witnesses $\mathrm{SOP}_1$ if there are
$\varepsilon>0$ and parameters $(a_\eta)_{\eta\in 2^{<\omega}}$ such that
$\{\varphi(x;a_{\sigma|n})=0:n<\omega\}$ is consistent for every
$\sigma\in2^\omega$, whereas
\[
 \inf_x\max\{\varphi(x;a_{\eta\frown1}),\varphi(x;a_\xi)\}
 \geqslant\varepsilon
 \qquad\text{whenever }\eta\frown0\unlhd\xi.
\]

\begin{fact}\label{fact:continuous-parallel}
Let $(a_i,b_i)_{i<\omega}$ be an indiscernible sequence in a continuous
theory, and suppose that
\[
 a_i\equiv_{a_{<i}b_{<i}}b_i\qquad(i<\omega).
\]
If $\varphi(x;y)$ is a continuous formula and $c\in\mathbb R$ are such
that $\{\varphi(x;a_i)\geqslant c:i<\omega\}$ is consistent and
$\{\varphi(x;b_i)\geqslant c:i<\omega\}$ is inconsistent, then the
theory has $\mathrm{SOP}_1$.
\end{fact}

\begin{proof}
This is the continuous version of the parallel sequences lemma
\cite[Lemma~2.28]{DK22}. We recall explicitly why that result
applies. Present a saturated continuous structure as a positive
structure by adding a relation $R_\psi$ for each zero condition
$\psi=0$. By \cite[Proposition~7.2.4 and Theorem~7.2.5]{Kamsma25},
the resulting positive theory is Hausdorff, hence thick, and its monster
has the same types and automorphisms,
and every positive formula is equivalent to a zero condition of a
continuous formula. Apply the parallel sequences lemma to
$R_{(c-\varphi)^+}$.
Here $t^+=\max\{t,0\}$.

It remains to check that the resulting positive $\mathrm{SOP}_1$
witness gives a continuous witness with a uniform gap. Write its
formula as $R_\theta(x;y)$, and let $\chi(y,z)$ be a positive
formula witnessing the inconsistency of the corresponding pairs.
Thus $\chi(y,z)$ is inconsistent with
$\exists x(R_\theta(x;y)\wedge R_\theta(x;z))$.
The latter formula is equivalent to the zero condition of
\[
 F(y,z)=\inf_x\max\{\theta(x;y),\theta(x;z)\}.
\]
By \cite[Proposition~7.2.8]{Kamsma25}, there is
$\varepsilon>0$ such that $\chi(y,z)$ implies $F(y,z)\geqslant
\varepsilon$. Thus the positive witness is also a continuous
$\mathrm{SOP}_1$ witness in the stated convention.
\end{proof}

Equip $\omega^{<\omega}$ with the language consisting of initial
segment, meet and lexicographic order, without level predicates.
We call $\eta_1,\ldots,\eta_k$ \emph{distant siblings} if, for some
$t\in\omega^{<\omega}$ and distinct $j_1,\ldots,j_k<\omega$, one has
$t\frown j_i\unlhd\eta_i$ for every $i$.

\begin{lemma}\label{lem:comb}
Let $(v_\eta)_{\eta\in\omega^{<\omega}}$ be strongly tree-indiscernible
in a continuous theory. Suppose that $P(x;y)$ is a continuous
formula, $c\in\mathbb R$ and $k\geqslant2$ satisfy the following:
\begin{enumerate}[(i)]
\item the conditions $P(x;v_\eta)\geqslant c$ are consistent along
every branch;
\item for any $k$ distant siblings $\eta_1,\ldots,\eta_k$, the
conditions $P(x;v_{\eta_j})\geqslant c$, $1\leqslant j\leqslant k$,
are inconsistent.
\end{enumerate}
Then the theory has $\mathrm{SOP}_1$.
\end{lemma}

\begin{proof}
For $m<\omega$, by an $m$-comb we mean a set
\[
 C=\{d_i:i<\omega\}\cup\{e_{i,j}:i<\omega,\ j<m\},
\]
where $d_0\lhd d_1\lhd\cdots$ and, for each $i<\omega$,
the nodes $e_{i,0},\ldots,e_{i,m-1},d_{i+1}$ lie in pairwise
distinct immediate-successor cones of $d_i$, with
\[
 e_{i,0}<_{\rm lex}\cdots<_{\rm lex}e_{i,m-1}
 <_{\rm lex}d_{i+1}.
\]
We call $(d_i)_{i<\omega}$ the spine and the nodes $e_{i,j}$
the teeth of the comb. For each $r<\omega$, the quantifier-free
tree type of $(d_i,e_{i,0},\ldots,e_{i,m-1})_{i<r}$
depends only on $m$ and $r$. Hence, by strong tree
indiscernibility and compactness, consistency of the conditions
indexed by an $m$-comb is independent of its choice.

We call a comb $C$ consistent if
$\{P(x;v_\eta)\geqslant c:\eta\in C\}$ is consistent.
A $0$-comb is consistent by (i), whereas a $(k-1)$-comb is
inconsistent by (ii): for any $i<\omega$, the nodes
$e_{i,0},\ldots,e_{i,k-2},d_{i+1}$ are $k$ distant siblings.
Choose $m\leqslant k-2$ maximal such that an $m$-comb is
consistent, and put $r=m+1$. For $i<\omega$ and $j<m$, define
\[
 t_i=(r,r)^i,\qquad a_i=t_i\frown r,\qquad b_i=t_i\frown m,
 \qquad s_{i,j}=t_i\frown j,
\]
where $(r,r)^i$ denotes the concatenation of $i$ copies of
$(r,r)$. Set
\[
 A_i=(v_{t_i},v_{s_{i,0}},\ldots,v_{s_{i,m-1}},v_{a_i}),\qquad
 B_i=(v_{t_i},v_{s_{i,0}},\ldots,v_{s_{i,m-1}},v_{b_i}).
\]
By strong tree indiscernibility, $(A_i,B_i)_{i<\omega}$ is
indiscernible and
\[
 A_i\equiv_{A_{<i}B_{<i}}B_i\qquad(i<\omega).
\]
Indeed, the quantifier-free tree type of any finite sequence
of the corresponding index tuples depends only on their order.
Moreover, replacing $a_i$ by $b_i$ preserves the quantifier-free
tree type over the indices occurring in $A_{<i}B_{<i}$, together
with $t_i,s_{i,0},\ldots,s_{i,m-1}$.

The set of indices occurring in $(A_i)_{i<\omega}$ is contained
in an $m$-comb with spine
\[
 t_0\lhd a_0\lhd t_1\lhd a_1\lhd\cdots
\]
and teeth $s_{i,j}$ at $t_i$ and $a_i\frown j$ at $a_i$,
for $j<m$. The set of indices occurring in $(B_i)_{i<\omega}$
is an $(m+1)$-comb with spine $(t_i)_{i<\omega}$ and teeth
$s_{i,0},\ldots,s_{i,m-1},b_i$ at $t_i$.
Thus the conditions indexed by the former set are consistent,
whereas those indexed by the latter are inconsistent, by the
choice of $m$.

For $w=(w_0,\ldots,w_{m+1})$, put
\[
 Q(x;w)=\min_{j<m+2}P(x;w_j).
\]
Then $\{Q(x;A_i)\geqslant c:i<\omega\}$ is consistent, while
$\{Q(x;B_i)\geqslant c:i<\omega\}$ is inconsistent.
Fact~\ref{fact:continuous-parallel} completes the proof.
\end{proof}

\section{\texorpdfstring{SOP$_1$}{SOP1} versus generic transitivity}
\label{sec:asymmetry}

\begin{lemma}\label{lem:separation}
Let $\nu$ be a global fim idempotent on a type-definable group $G$, and
suppose that $\nu\ne\check\nu$. There are a small model $M$ over which $\nu$ is fim
and an $M$-definable finite-valued function $h:G\to[0,1]$
with rational values
such that
\begin{equation}\label{eq:gap}
 A:=\nu(h)>B:=\sup_{z\in G(\FC)}\nu_y\bigl(h(zy^{-1})\bigr).
\end{equation}
\end{lemma}

\begin{proof}
Put $\mu=\nu*\check\nu$. By \eqref{eq:inverse-convolution},
$\check\mu=\mu$, so $\nu\ne\mu$.
Choose a relatively definable set $D\subseteq G$ such that
$\nu(D)>\mu(D)$.
Choose a small model $M$ containing the parameters of $D$ and
over which $\nu$ is fim.

On the real space $C(S_G(M))$, define
\[
 (Rf)(z)=\nu_y\bigl(f(zy^{-1})\bigr).
\]
For locally constant functions the output is continuous by
definability of $\nu$. Uniform approximation extends $R$ to a linear operator on
$C(S_G(M))$ satisfying $R1=1$, $Rf\geqslant0$ whenever
$f\geqslant0$, and $\|Rf\|_\infty\leqslant\|f\|_\infty$.
Idempotence of $\nu$ and symmetry of $\nu_x\otimes\nu_y$ give
\[
 R^2=R,\qquad \nu(R\mathbf{1}_D)=\mu(\mathbf{1}_D).
\]
For example, the two successive averages in $R^2f(z)$ evaluate
$f(za^{-1}b^{-1})=f(z(ba)^{-1})$, and the product $ba$ has law $\nu$.
The second equality uses symmetry of $\nu_x\otimes\nu_y$
to interchange the two integrations.

The function
\[
 f=\frac{1+\mathbf{1}_D-R\mathbf{1}_D}{2}
\]
takes values in $[0,1]$ and satisfies
\[
 Rf=\frac12,\qquad
 \nu(f)=\frac12+\frac{\nu(\mathbf{1}_D)-\mu(\mathbf{1}_D)}2>\frac12.
\]
Choose $\delta>0$ with $2\delta<\nu(f)-\tfrac12$.
Since rational-valued locally constant functions are uniformly dense
in $C(S_G(M))$, there is such an $h$ taking values in $[0,1]$ with
$\|h-f\|_\infty<\delta$. It follows that
\[
 \nu(h)>\nu(f)-\delta>\frac12+\delta
 \geqslant\sup_{z\in G(\FC)}Rh(z),
\]
which proves \eqref{eq:gap}.
\end{proof}

\begin{proposition}\label{prop:asymmetry-tree}
If a type-definable group in $T$ carries a global fim idempotent $\nu$
such that $\nu\ne\check\nu$, then $T^R$ has $\SOP$.
\end{proposition}

\begin{proof}
Choose $M,h,A,B$ as in Lemma~\ref{lem:separation}, and write
$G=\pi(\FC)$, where $\pi$ is a partial type over $M$.
Work in a sufficiently saturated model of $T^R$ containing the
constant copy $M^c$ of $M$. Throughout the proof, expressions
involving $h$ and the group operations are interpreted using
their fixed ambient definable extensions.
Let $I=\{*\}\sqcup\omega^{<\omega}$.

We first choose a family $(L_i)_{i\in I}$ such that, for every
$n\geqslant1$, every ordered tuple of pairwise distinct indices
$i_1,\ldots,i_n\in I$, and every formula
$\varphi(x_1,\ldots,x_n;m)\in\CL(M)$, we have
\[
 \E[\varphi(L_{i_1},\ldots,L_{i_n};m)]
 =
 \nu^{(n)}\bigl(\varphi(x_1,\ldots,x_n;m)\bigr).
\]
To verify consistency of these requirements, note that
\[
 (\operatorname{pr}_{j_1,\ldots,j_r})_*
       (\nu^{(n)}|M)=\nu^{(r)}|M
\]
whenever $j_1,\ldots,j_r$ are pairwise distinct elements of
$\{1,\ldots,n\}$.
Indeed, the Morley powers are symmetric by
\cite[Theorems~2.18 and~5.16(a)]{InTheWild}.
Consequently, for any finite set of indices, the randomized
type corresponding to the appropriate Morley power satisfies
all the requirements involving those indices.
Compactness and saturation therefore give the desired family.
Since each $L_i$ has law $\nu|M$ and $\nu$ concentrates on $G$,
we have $L_i\in G^R$ for every $i\in I$.

Define
\begin{equation}\label{eq:tree-products}
 V_\varnothing=L_\varnothing,
 \qquad V_{\eta\frown i}=L_{\eta\frown i}V_\eta.
\end{equation}
Thus, if $d=|\eta|$, then
\[
 V_\eta=L_{\eta|d}L_{\eta|(d-1)}\cdots L_{\eta|0},
\]
where $\eta|0=\varnothing$.

\begin{clm}
The tree $(V_\eta)_{\eta\in\omega^{<\omega}}$ is strongly
tree-indiscernible over $M^c$.
\end{clm}

\begin{clmproof}
Fix a finite tuple $\bar\eta=(\eta_1,\ldots,\eta_n)$, let $F$
be its closure under meets.
Let $t_0$ be the least element of $F$.
For each $t\in F\setminus\{t_0\}$, let $p(t)$ be the greatest
proper predecessor of $t$ belonging to $F$.
Define
\[
 \mathcal B_{t_0}:=\{\rho:\rho\unlhd t_0\},
 \qquad
 \mathcal B_t:=\{\rho:p(t)\lhd\rho\unlhd t\}
 \quad(t\in F\setminus\{t_0\}),
\]
and
\[
 Z_t:=L_{\rho_m}\cdots L_{\rho_0}
 \qquad
 \text{where }\mathcal B_t=\{\rho_0\lhd\cdots\lhd\rho_m\}.
\]
The blocks $\mathcal B_t$ are nonempty and pairwise disjoint,
since $F$ is closed under meets.
By \eqref{eq:tree-products}, for every $t\in F$,
\[
 V_t=Z_{s_\ell}\cdots Z_{s_0}
 \qquad
 \text{where }
 \{s\in F:s\unlhd t\}
 =\{s_0\lhd\cdots\lhd s_\ell=t\}.
\]

Now let $\bar\eta'=(\eta'_1,\ldots,\eta'_n)$ have the same
quantifier-free tree type as $\bar\eta$.
There is an induced isomorphism
\[
 f:F\longrightarrow F',
 \qquad f(\eta_i)=\eta'_i,
\]
where $F'$ is the meet closure of $\bar\eta'$.
Apply the same construction to $F'$, using primed notation.
For any enumeration $F=\{t_1,\ldots,t_r\}$,
Lemma~\ref{lem:blocks} gives
\[
 \E[\theta(Z_{t_1},\ldots,Z_{t_r})]
 =
 \nu^{(r)}(\theta)
 =
 \E[\theta(Z'_{f(t_1)},\ldots,Z'_{f(t_r)})]
 \qquad
 \bigl(\theta(z_1,\ldots,z_r)\in\CL(M)\bigr).
\]
By quantifier elimination for randomization
\cite[Theorem~2.9]{BenYaacov2009}, the two tuples of block
products have the same type over $M^c$.
Since $f$ preserves the initial-segment order, the reconstruction
formula for $V_t$ gives
\[
 (V_{\eta_1},\ldots,V_{\eta_n})
 \equiv_{M^c}
 (V_{\eta'_1},\ldots,V_{\eta'_n}).
\]
This proves the claim.
\end{clmproof}

Put
\[
 P(X;V):=\E[h(XV^{-1})],
 \qquad \varepsilon:=\frac{A-B}{8},
 \qquad c:=\frac{A+B}{2}=B+4\varepsilon.
\]
Since $h$ is finite-valued and $M$-definable, $P$ is a
$[0,1]$-valued continuous formula over $M^c$.

\begin{clm}\label{clm:asymmetry-branches}
For every $\sigma\in\omega^\omega$, the partial type
\[
 \pi^R(X)\cup
 \{P(X;V_{\sigma|n})=A:n<\omega\}
\]
is consistent.
\end{clm}

\begin{clmproof}
Fix $\tau\in\omega^{<\omega}$, put $d=|\tau|$, and let
$X=L_*V_\tau\in G^R$.
For every $\eta\unlhd\tau$, we have
\[
 L_*V_\tau V_\eta^{-1}
 =
 \begin{cases}
 L_*L_{\tau|d}\cdots L_{\tau|(|\eta|+1)},
     &\eta\lhd\tau,\\
 L_*,&\eta=\tau.
 \end{cases}
\]
Each expression on the right is a nonempty product of labels
with pairwise distinct indices. Lemma~\ref{lem:blocks} gives
\begin{equation}\label{eq:branch-value}
 P(X;V_\eta)
 =\E[h(L_*V_\tau V_\eta^{-1})]
 =\nu(h)=A
 \qquad(\eta\unlhd\tau).
\end{equation}
The claim follows by compactness.
\end{clmproof}

\begin{clm}\label{clm:asymmetry-forks}
There are $k\geqslant2$ and a finite conjunction $\chi$ of
formulas from $\pi$ such that, for every ambient randomized
tuple $X$ and every $k$ distant siblings $\eta_1,\ldots,\eta_k$,
\[
 \E[\chi(X)]=1
 \quad\Longrightarrow\quad
 \frac1k\sum_{i=1}^kP(X;V_{\eta_i})\leqslant B+2\varepsilon.
\]
\end{clm}

\begin{clmproof}
Choose $N\geqslant1$ such that $Nh$ is integer-valued, and put
\[
 \psi_j(z):=(h(z)\geqslant j/N),
 \qquad
 \varphi_j(y;x,v):=\psi_j\bigl(x(yv)^{-1}\bigr)
 \quad(1\leqslant j\leqslant N).
\]
Thus each $\psi_j$ is an $M$-formula and
\[
 h(z)=\frac1N\sum_{j=1}^N\mathbf1_{\psi_j(z)}.
\]
Apply Definition~\ref{def:fim} to $\varphi_1,\ldots,\varphi_N$.
For a common sufficiently large $k\geqslant2$, there are
$M$-definable sets $\Theta_j\subseteq(\FC^x)^k$ such that
\[
 \nu^{(k)}(\Theta_j)>1-\frac{\varepsilon}{N}
\]
and
\[
 \left|
 \frac1k\sum_{i=1}^k\mathbf1_{\varphi_j(y_i;a,b)}
 -\nu_y(\varphi_j(y;a,b))
 \right|<\varepsilon
 \qquad(a,b\in\FC^x,\ \bar y\in\Theta_j).
\]
Put $\Theta:=\bigcap_{j=1}^N\Theta_j$. Then $\nu^{(k)}(\Theta)>1-\varepsilon$,
and averaging the preceding estimates gives
\[
 \left|
 \frac1k\sum_{i=1}^k h\bigl(a(y_i b)^{-1}\bigr)
 -\nu_y\bigl(h(a(yb)^{-1})\bigr)
 \right|<\varepsilon
 \qquad(a,b\in\FC^x,\ \bar y\in\Theta).
\]
For $a,b\in G$, concentration of $\nu$ on $G$ gives
\[
 \nu_y\bigl(h(a(yb)^{-1})\bigr)
 =\nu_y\bigl(h((ab^{-1})y^{-1})\bigr)
 \leqslant B.
\]
Consequently,
\[
 \pi(x)\cup\pi(v)\cup\{\Theta(\bar y)\}
 \models
 \frac1k\sum_{i=1}^k h\bigl(x(y_i v)^{-1}\bigr)
 \leqslant B+\varepsilon.
\]
The inequality is a first-order condition involving finitely
many coordinates, since $h$ has finite range.
By compactness, taking the conjunction of the finite fragments
of $\pi$ obtained in the $x$- and $v$-variables, we obtain a
finite conjunction $\chi$ of formulas from $\pi$ such that
\[
 \chi(x)\wedge\chi(v)\wedge\Theta(\bar y)
 \quad\Longrightarrow\quad
 \frac1k\sum_{i=1}^k h\bigl(x(y_i v)^{-1}\bigr)
 \leqslant B+\varepsilon.
\]

Let $\eta_1,\ldots,\eta_k$ be distant siblings with common
meet $t$, and put
\[
 Y_i:=L_{\eta_i|d_i}\cdots L_{\eta_i|(|t|+1)},
 \qquad d_i:=|\eta_i|
 \quad(1\leqslant i\leqslant k).
\]
The sets $\{\rho:t\lhd\rho\unlhd\eta_i\}$ are nonempty and
pairwise disjoint. Hence \eqref{eq:tree-products} and
Lemma~\ref{lem:blocks} give
\[
 V_{\eta_i}=Y_iV_t,
 \qquad
 \E[\Theta(Y_1,\ldots,Y_k)]=\nu^{(k)}(\Theta).
\]
Fix $X$ with $\E[\chi(X)]=1$.
Since $V_t\in G^R$, we also have $\E[\chi(V_t)]=1$.
Using the preceding first-order implication in the
randomization \cite[Section~2]{BenYaacov2009} and
$0\leqslant h\leqslant1$, we obtain
\begin{align}
 \frac1k\sum_{i=1}^kP(X;V_{\eta_i})
 &=\E\left[
      \frac1k\sum_{i=1}^k h\bigl(X(Y_iV_t)^{-1}\bigr)
     \right]\notag\\
 &\leqslant B+\varepsilon+1
   -\E[\chi(X)\wedge\chi(V_t)\wedge
        \Theta(Y_1,\ldots,Y_k)]\notag\\
 &=B+\varepsilon+1-\nu^{(k)}(\Theta)
 <B+2\varepsilon.
 \label{eq:fork-bound}
\end{align}
This proves the claim.
\end{clmproof}

Choose $k$ and $\chi$ as in Claim~\ref{clm:asymmetry-forks},
and define
\[
 Q(X;V):=\max\{1-\E[\chi(X)],(c-P(X;V))^+\}.
\]
Since $A>c>B+2\varepsilon$, Claims~\ref{clm:asymmetry-branches}
and~\ref{clm:asymmetry-forks} imply that
\[
 \{Q(X;V_{\sigma|n})=0:n<\omega\}
 \quad\text{is consistent for every }\sigma\in\omega^\omega,
\]
whereas
\[
 \{Q(X;V_{\eta_i})=0:1\leqslant i\leqslant k\}
 \quad\text{is inconsistent for every }k\text{ distant siblings.}
\]

Write
\[
 Q(X;V)=Q_0(X_0;V_0,a),
\]
where $X_0,V_0$ are finite subtuples and $a$ is a finite tuple
from $M^c$. Put $v_\eta:=(V_{\eta,0},a)$, where
$V_{\eta,0}$ denotes the corresponding subtuple of $V_\eta$.
By Claim~1, $(v_\eta)_{\eta\in\omega^{<\omega}}$ is strongly
tree-indiscernible. The preceding consistency and inconsistency
conditions hold for $Q_0(X_0;v_\eta)$, since every assignment
to $X_0$ extends to the unused ambient coordinates of $X$.
Lemma~\ref{lem:comb}, applied to $1-Q_0$ at threshold $1$,
yields $\SOP$ in $T^R$.
\end{proof}

\begin{proof}[Proof of Theorem~\ref{thm:main}]
Suppose that $T^R$ is $\NSOP$. By
Proposition~\ref{prop:asymmetry-tree}, every global fim idempotent
on $G$ is invariant under inversion. Proposition~\ref{prop:symmetry-reduction}
then implies that every such measure is generically transitive.
The converse and the remaining conclusions are
Fact~\ref{fact:classification}.
\end{proof}

\begin{proof}[Proof of Corollary~\ref{cor:ground-main}]
By Theorem \ref{thm:randomization-preserves-nsop1} and Theorem~\ref{thm:main}.
\end{proof}

For a global fim idempotent $\mu$, write $S_\mu=\supp(\mu)$.
The product of invariant types makes $S_\mu$ a compact Hausdorff
semigroup with continuous right translations
\cite[Proposition~4.3 and Corollary~4.4]{Artem_Kyle}.
It is \emph{left-simple} if $S_\mu*p=S_\mu$ for every $p\in S_\mu$.

\begin{cor}\label{cor:support}
Suppose that $T^R$ is $\NSOP$. For every type-definable group $G$ in $T$
and every global fim idempotent $\mu$ on $G$, one has
\[
 \mu*p=\mu\quad\text{and}\quad S_\mu*p=S_\mu
 \qquad(p\in S_\mu).
\]
In particular, $S_\mu$ is its unique minimal left ideal.
The same conclusion holds for $T$ being $\NSOP$.
\end{cor}
\begin{proof}
By Theorem~\ref{thm:main}, the measure $\mu$ concentrates on its
right stabilizer $H$. Fix $p\in S_\mu$ and a definable set $D$,
and choose a small model $N$ containing the parameters of $D$
and over which $\mu$ is fim. If $b\models p|N$, then $b\in H$, and hence
\[
 (\mu*p)(D)=(\mu*\delta_b)(D)=\mu(D).
\]
We have proved $\mu*p=\mu$ for every $p\in S_\mu$.
By \cite[Proposition~3.51]{CGK}, this is equivalent to
$S_\mu*p=S_\mu$ for every $p\in S_\mu$. Finally, any nonempty left ideal
containing $p$ contains $S_\mu*p=S_\mu$, proving the assertion
about minimal left ideals.
\end{proof}

For types, generic transitivity of generically stable idempotents
was established in simple theories and in rosy theories in
\cite[Propositions~2.34 and~2.37]{CGK}, respectively.
The following corollary records the consequence of our measure
theorem for types.

\begin{cor}\label{cor:types}
Suppose that $T^R$ is $\NSOP$, and let $p$ be a global fim
idempotent type on a type-definable group $G$. Then $p$ is generically
transitive. Explicitly, if $M$ is a small model containing the group parameters
and $p$ is fim over $M$, then
\[
 (a,b)\models p^{(2)}|M
 \quad\Longrightarrow\quad (ab,b)\models p^{(2)}|M.
\]
Moreover, $p$ is the unique left-invariant and the unique
right-invariant global type concentrating on its stabilizer.
\end{cor}
\begin{proof}
Apply Theorem~\ref{thm:main} to the Dirac measure $\delta_p$ and
use \eqref{eq:shear}.
\end{proof}

\begin{cor}\label{cor:failure}
If some type-definable group in $T$ carries a global fim idempotent
which is not generically transitive, then $T$ has $\SOP$ and IP.
\end{cor}
\begin{proof}
The contrapositive of Proposition~\ref{prop:symmetry-reduction}
produces a global fim idempotent on the same group which is not
invariant under inversion. Proposition~\ref{prop:asymmetry-tree}
and Theorem \ref{thm:randomization-preserves-nsop1} apply and we obtain $\SOP$.
The IP conclusion follows by \cite{ChernikovNIP2026}.
\end{proof}

\appendix
\section{Randomization preserves \texorpdfstring{$\mathrm{NSOP}_1$}{NSOP1}}
\label{app:preservation}

\providecommand{\llbracket}{[\![}
\providecommand{\rrbracket}{]\!]}

We give a proof of the following theorem of Ben Yaacov, Chernikov
and Ramsey, announced in \cite{BenYaacovNSOP2019}
(note that \cite{BYCR2024} is listed as a 2024 preprint
in the bibliography of \cite{CGK}).
A related combination of
pointwise independence and independence in the event algebra
was studied, for stable theories, in
\cite[Section~7.3]{AndrewsGoldbringKeisler2019}.

\begin{theorem}[Ben Yaacov--Chernikov--Ramsey]
\label{thm:randomization-preserves-nsop1}
Let $T$ be a complete first-order theory, and let $T^R$
be its Keisler randomization. Then
\[
 T\text{ is }\NSOP
 \quad\Longleftrightarrow\quad
 T^R\text{ is }\NSOP.
\]
\end{theorem}

We use Kim's lemma, symmetry and the independence theorem for
Kim-independence over models of an $\NSOP$ theory
\cite{KaplanRamsey2020}. We write $a\ind^K_M b$ for
Kim-independence. We use disintegration on standard Borel
spaces \cite[Theorem~8.5]{Kallenberg2021}, Souslin's theorem
\cite[Theorem~14.11]{Kechris1995}, and Jankov--von Neumann
uniformization in its universally measurable form
\cite[Theorem~29.9]{Kechris1995}.
In continuous structures, finite satisfiability is understood
approximately: every finite family of conditions can be
approximated arbitrarily closely by a tuple from the base model.

We first prove a criterion for $\NSOP$ in terms of an auxiliary
independence relation over separable models. For countable $T$,
we then define such a relation using Borel Kim-independence
and conditional independence of event algebras.
Lemma~\ref{lem:app-random-coheir} shows that coheir implies
this relation, and Proposition~\ref{prop:app-random-amalgamation}
establishes independent amalgamation.
The final subsection removes the countability
assumption.

\renewcommand{\theHclm}{\theHtheorem.\arabic{clm}}

\subsection{A separable coheir criterion}

\begin{lemma}\label{lem:app-continuous-criterion}
Let $U$ be a complete continuous theory in a countable language.
Suppose that an invariant relation $\ind^*$ is defined for finite
tuples over separable elementary models and satisfies:
\begin{enumerate}[(i)]
\item if $\tp(a/Nb)$ is finitely satisfiable in $N$, then
$a\ind^*_N b$;
\item $a\ind^*_N b$ implies $b\ind^*_N a$;
\item if $c_1\ind^*_N c_2$,
$b_1\ind^*_N c_1$, $b_2\ind^*_N c_2$, and
$b_1\equiv_Nb_2$, there is $b$ with
$b\equiv_{Nc_1}b_1$ and $b\equiv_{Nc_2}b_2$.
\end{enumerate}
Then $U$ is $\NSOP$.
\end{lemma}

\begin{proof}
Suppose that $U$ has $\SOP$, witnessed by $\psi(x;y)$,
with consistent zero conditions along branches and a rational
inconsistency gap $\eta>0$. Put
\[
 F(y,z):=\inf_x\max\{\psi(x;y),\psi(x;z)\},
\]
and take rational $0<\delta<\eta$. 
We will apply \cite[Proposition~A.7]{HK21} to an ordinary
first-order presentation of $U$ and transfer the resulting
configuration to a separable elementary model of $U$.

Fix a monster $\mathcal U\models U$. Choose a countable family
$\mathcal F$ of continuous formulas, uniformly dense in all
formulas, containing $\psi$ and all atomic formulas, and closed
under variable substitution, $\min$, $\max$, complementation,
rational truncated subtraction and $\inf$-quantification.
Expand the underlying sets and functions of $\mathcal U$ by
relations
\[
 R_{\theta\leq r}(a)
 \quad\Longleftrightarrow\quad
 \theta^{\mathcal U}(a)\leq r
 \qquad
 (\theta\in\mathcal F,\ r\in\mathbb Q\cap[0,1]),
\]
and let $U_\#$ be the complete ordinary first-order theory of
this countable presentation. Write
$R_{\theta\geq r}:=R_{1-\theta\leq1-r}$.

\begin{clm}\label{clm:app-reconstruction}
For $A\models U_\#$, the values
\[
 \theta_A(a):=
 \inf\{r\in\mathbb Q\cap[0,1]:
                  A\models R_{\theta\leq r}(a)\}
\]
determine a continuous model $\widehat A\models U$ after
quotienting by distance zero and completing. Moreover,
\[
 A_0\preccurlyeq A
 \quad\Longrightarrow\quad
 \widehat{A_0}\preccurlyeq\widehat A.
\]
\end{clm}

\begin{clmproof}
The rational comparison sentences of $U_\#$ give
\[
 \theta_A(a)<r
 \quad\Longrightarrow\quad
 A\models R_{\theta\leq r}(a)
 \quad\Longrightarrow\quad
 \theta_A(a)\leq r,
\]
and, for each connective occurring in $\mathcal F$,
\[
 \bigl(g(\theta_1,\ldots,\theta_k)\bigr)_A(a)
 =g\bigl((\theta_1)_A(a),\ldots,(\theta_k)_A(a)\bigr).
\]
For quantifiers, the sentences
\[
\begin{aligned}
 R_{\inf_y\theta\leq r}(x)
 &\ \longrightarrow\ \exists y\,R_{\theta\leq s}(x,y)
 &&(r<s),\\
 R_{\theta\leq r}(x,y)
 &\ \longrightarrow\ R_{\inf_y\theta\leq r}(x)
\end{aligned}
\]
give
\[
 (\inf_y\theta)_A(a)=\inf_{b\in A}\theta_A(a,b).
\]
The same rational comparisons give the pseudometric axioms
and the prescribed continuity moduli of all symbols. They also
give $\sigma_A=\sigma^{\mathcal U}$ for every sentence
$\sigma\in\mathcal F$. Thus the quotient completion satisfies
$U$, by uniform approximation using $\mathcal F$.
If $A_0\preccurlyeq A$, then
\[
 \theta_{A_0}(a)=\theta_A(a)
 \qquad(a\in A_0,\ \theta\in\mathcal F).
\]
Uniform approximation and continuity extend this equality to
all formulas and all tuples from $\widehat{A_0}$, proving
$\widehat{A_0}\preccurlyeq\widehat A$.
\end{clmproof}

\begin{clm}
There are a separable elementary model $N\models U$ and finite
tuples $b_1,b_2,c_1,c_2$ such that
\[
\begin{gathered}
 \tp(c_1/Nc_2),\quad \tp(c_1/Nb_1),\quad
 \tp(c_2/Nb_2)\text{ are finitely satisfiable in }N,\\
 b_1\equiv_Nb_2,\qquad
 \psi(b_i;c_i)\leq\delta\quad(i=1,2),\qquad
 F(c_1,c_2)\geq\eta.
\end{gathered}
\]
\end{clm}

\begin{clmproof}
In $U_\#$, put
\[
 \Phi(x;y):=R_{\psi\leq\delta}(x;y),
 \qquad
 \Psi(y,z):=R_{F\geq\eta}(y,z).
\]
The given tree witnesses $\SOP$ for $\Phi$, and
\[
 U_\#\models
 \neg\exists x,y,z\,
 \bigl(\Phi(x;y)\wedge\Phi(x;z)\wedge\Psi(y,z)\bigr).
\]
Take the full ordinary first-order Morleyization of $U_\#$.
This is a countable definitional expansion with quantifier
elimination; it is inductive and all its models are
existentially closed. Apply \cite[Proposition~A.7]{HK21} and
take reducts to $U_\#$. We obtain $N_0\preccurlyeq A\models U_\#$
and tuples satisfying
\[
\begin{gathered}
 \tp_{U_\#}(c_1/N_0c_2),\quad
 \tp_{U_\#}(c_1/N_0b_1),\quad
 \tp_{U_\#}(c_2/N_0b_2)
 \text{ are finitely satisfiable in }N_0,\\
 b_1\equiv^{U_\#}_{N_0}b_2,\qquad
 A\models\Phi(b_1;c_1)\wedge\Phi(b_2;c_2)\wedge\Psi(c_1,c_2).
\end{gathered}
\]
Here $N_0$ is countable: the proof takes the Skolem hull of
the entries of an $\omega$-sequence in a countable expansion.
It also retains the specified predicate $\Psi$.
By Claim~\ref{clm:app-reconstruction},
\[
 N:=\widehat{N_0}\preccurlyeq\widehat A,
 \qquad N\text{ is separable}.
\]
We suppress the quotient notation on tuples. Equality of
first-order types and continuity give
\[
 b_1\equiv^{U_\#}_{N_0}b_2
 \quad\Longrightarrow\quad
 b_1\equiv^U_Nb_2.
\]

To transfer finite satisfiability, suppose that
$\tp_{U_\#}(a/N_0c)$ is finitely satisfiable in $N_0$.
Fix continuous formulas $\theta_j(u;c,n)$, where $n\in N$,
and $\varepsilon>0$. Choose $n_0\in N_0$ and
$\theta_j'\in\mathcal F$ such that
\[
 \sup_u\bigl|
 \theta_j(u;c,n)-\theta_j'(u;c,n_0)
 \bigr|<\varepsilon/3
 \qquad(1\leqslant j\leqslant k).
\]
Choose $r_j,s_j\in\mathbb Q\cap[0,1]$ with
$s_j-r_j<\varepsilon/3$ such that
\[
\begin{aligned}
 r_j=0&\quad\text{or}\quad r_j<\theta_j'(a;c,n_0),\\
 s_j=1&\quad\text{or}\quad \theta_j'(a;c,n_0)<s_j.
\end{aligned}
\]
Then
\[
 A\models\bigwedge_{j=1}^k
 \bigl(R_{\theta_j'\geq r_j}(a;c,n_0)
       \wedge R_{\theta_j'\leq s_j}(a;c,n_0)\bigr).
\]
Finite satisfiability gives $a_0\in N_0$ satisfying this
conjunction. Consequently,
\[
\begin{aligned}
 \bigl|\theta_j(a_0;c,n)-\theta_j(a;c,n)\bigr|
 &<\frac{2\varepsilon}{3}
   +\bigl|\theta_j'(a_0;c,n_0)-\theta_j'(a;c,n_0)\bigr|\\
 &\leqslant\frac{2\varepsilon}{3}+s_j-r_j
 <\varepsilon.
\end{aligned}
\]
Thus $\tp_U(a/Nc)$ is approximately finitely satisfiable in
$N$. Applying this to the three displayed types gives the
required finite satisfiability. Finally,
\[
 \psi(b_i;c_i)\leq\delta\quad(i=1,2),
 \qquad F(c_1,c_2)\geq\eta
\]
follow directly from the defining comparison relations.
\end{clmproof}

By (i),
\[
 c_1\ind^*_N c_2,
 \qquad c_1\ind^*_N b_1,
 \qquad c_2\ind^*_N b_2.
\]
By (ii), $b_i\ind^*_N c_i$ for $i=1,2$.
Since $b_1\equiv_Nb_2$, condition (iii) gives $b$ such that
\[
 b\equiv_{Nc_i}b_i
 \quad\Longrightarrow\quad
 \psi(b;c_i)=\psi(b_i;c_i)\leq\delta
 \qquad(i=1,2).
\]
Hence
\[
 \eta\leq F(c_1,c_2)
 \leq\max\{\psi(b;c_1),\psi(b;c_2)\}
 \leq\delta<\eta,
\]
a contradiction.
\end{proof}

\subsection{Enumerations of models and Borel Kim-independence}

We allow many-sorted classical theories. Until the final proof,
assume that $T$ is countable and $\NSOP$, and fix a sort
$\mathsf S$ such that $T\models\exists u,v\in\mathsf S\ (u\ne v)$.
Write $\mathsf S^R$ for the corresponding random-element sort.
The Stone spaces of $T$ in countable tuples are compact
metrizable. 
Let $X=(X_i)_{i<\omega}$ be a countable tuple of random elements,
and let $x=(x_i)_{i<\omega}$ be the corresponding tuple of
classical variables.
For a classical formula $\theta(x)$, write
$\llbracket\theta(X)\rrbracket$ for the event that $\theta$
holds of $X$. Its probability is
$\mathbb P\llbracket\theta(X)\rrbracket=\E[\theta(X)]$.
Define the Borel probability measure $\mu_X$
on $S_x(T)$ by
\[
 \mu_X([\theta(x)])
 =\mathbb P\llbracket\theta(X)\rrbracket
 \qquad(\theta(x)\in\CL).
\]
Conversely, every Borel probability measure on $S_x(T)$ determines
a randomized type
\cite[Theorem~2.9 and Corollary~2.10]{BenYaacov2009}.
For a Borel set $B\subseteq S_x(T)$, we say that the classical
type of $X$ belongs to $B$ almost surely if $\mu_X(B)=1$.

\begin{lemma}\label{lem:app-model-code}
Let $\mathcal N\models T^R$ be separable. There is a countable
tuple $\bar M=(M_i)_{i<\omega}$ of its random elements such that:
\begin{enumerate}[(i)]
\item its coordinates in each random-element sort are dense in
that sort of $\mathcal N$;
\item almost surely, its range, taken sort by sort, is an
elementary model of $T$;
\item the events $\llbracket\theta(\bar M)\rrbracket$ generate
the completed event algebra $\mathcal B_{\mathcal N}$.
\end{enumerate}
Let $\bar m=(m_i)_{i<\omega}$ be the corresponding tuple of
classical variables. The set $\mathcal M\subseteq S_{\bar m}(T)$
of types of enumerations of elementary submodels is Borel.
\end{lemma}

\begin{proof}
Starting with a countable tuple whose coordinates are dense in
each random-element sort, repeatedly add, for each
chosen finite tuple $\bar d$ and each formula $\varphi(u;\bar d)$,
a witness $c\in\mathcal N$ satisfying
\[
 \llbracket\exists u\,\varphi(u;\bar d)\rrbracket
 =\llbracket\varphi(c;\bar d)\rrbracket.
\]
Such witnesses exist by fullness and completeness
\cite[Theorem~2.7]{BenYaacov2009}.
There are countably many requirements, so the Tarski--Vaught
test gives (ii) outside a countable union of null sets.
Fullness and the choice of $\mathsf S$ imply that every event
in $\mathcal B_{\mathcal N}$ has the form
$\llbracket u=v\rrbracket$ for some
$u,v\in\mathsf S^R(\mathcal N)$.
Approximating $u,v$ by coordinates of $\bar M$ in this sort
gives (iii).

Finally,
\[
 \mathcal M=
 \bigcap_{\theta,\bar i}
 \left(
 [\neg\exists u\,\theta(u;\bar m_{\bar i})]
 \cup\bigcup_{\bar j}[\theta(\bar m_{\bar j};\bar m_{\bar i})]
 \right),
\]
where $\bar j$ ranges over tuples of indices of the appropriate
sorts. Thus $\mathcal M$ is a countable intersection of open sets.
\end{proof}

Put $\kappa=\mu_{\bar M}$, so $\kappa(\mathcal M)=1$.
For $z\in\mathcal M$, let $M_z$ be the model enumerated by a
realization of $z$, up to the indicated enumeration.

\begin{lemma}\label{lem:app-borel-kim}
For each formula $\varphi(x;y,\bar m_0)$, the locus
\[
 \{\tp(\bar m,b):z=\tp(\bar m)\in\mathcal M,
       \ \varphi(x;b,\bar m_0)\text{ Kim-divides over }M_z\}
\]
is Borel. Consequently, for countable tuples $a,b$, the condition
$a\ind^K_{M_z} b$ is Borel on the space of types
$\tp(\bar m,a,b)$ with $z\in\mathcal M$.
\end{lemma}

\begin{proof}
Let $\mathcal H$ be the set of types $\tp(\bar m,(b_i)_{i<\omega})$
such that $z=\tp(\bar m)\in\mathcal M$ and
\[
 \begin{gathered}
 (b_i)_{i<\omega}\text{ is indiscernible over }M_z,\\
 \tp(b_i/M_zb_{<i})\text{ is finitely satisfiable in }M_z
 \qquad(i<\omega).
 \end{gathered}
\]
Indiscernibility is closed. Finite satisfiability is the
intersection, over $i,\theta,\bar k$, of the open sets
\[
 [\neg\theta(b_i;\bar m_{\bar k},b_{<i})]
 \cup\bigcup_{\bar j}
 [\theta(\bar m_{\bar j};\bar m_{\bar k},b_{<i})].
\]
Hence $\mathcal H$ is Borel.

Every sequence in $\mathcal H$ is a Morley sequence in a single
global coheir. Indeed, in common variables $y$, put
\[
 \Sigma(y)=\bigcup_{i<\omega}\tp(b_i/M_zb_{<i}),
\]
replacing only the distinguished tuple $b_i$ by $y$.
For finite $\Delta\subseteq\Sigma$, choose $n$ such that
\[
 \Delta\subseteq\bigcup_{i<n}\tp(b_i/M_zb_{<i}).
\]
Indiscernibility gives $b_n\models\Delta$, and finite
satisfiability of $\tp(b_n/M_zb_{<n})$ gives a realization in
$M_z$. Thus $\Sigma$ extends to a global coheir $q$, with
\[
 b_i\models q|M_zb_{<i}\qquad(i<\omega).
\]

Let $\pi$ be restriction to $(\bar m,b_0)$ and put
\[
 \begin{gathered}
 \mathcal X=\{\tp(\bar m,b):\tp(\bar m)\in\mathcal M\},\\
 \mathcal C_\varphi=
 \mathcal H\cap\bigcap_{n\geqslant1}
 \left[\exists x\bigwedge_{i<n}\varphi(x;b_i,\bar m_0)\right].
 \end{gathered}
\]
Every type over a model has a global coheir extension, so
$\pi(\mathcal H)=\mathcal X$. By Kim's lemma, the Kim-dividing locus
$D_\varphi\subseteq\mathcal X$ satisfies
\[
 D_\varphi=\pi(\mathcal H\setminus\mathcal C_\varphi),
 \qquad \mathcal X\setminus D_\varphi=\pi(\mathcal C_\varphi).
\]
Both sets are analytic, hence Borel by Souslin's theorem
\cite[Theorem~14.11]{Kechris1995}.
Finally, Kim-forking equals Kim-dividing over models
\cite[Proposition~3.19]{KaplanRamsey2020}; taking a countable
intersection over formulas and finite subtuples gives the
assertion about Kim-independence.
\end{proof}

\subsection{The auxiliary independence relation}

For a countable random tuple $X$, let $\mathcal B_{X\mathcal N}$
be the completed probability algebra generated by
$\mathcal B_{\mathcal N}$ and
$\{\llbracket\varphi(X,\bar M)\rrbracket:\varphi\in\CL\}$.
For countable random tuples $X,Y$, write $x,y$ for the
corresponding tuples of classical variables. Define
\[
 X\ind^*_{\mathcal N} Y
 \quad\Longleftrightarrow\quad
 \begin{cases}
 \mathbb E[\mathbf1_E\mathbf1_F\mid\mathcal B_{\mathcal N}]
 =\mathbb E[\mathbf1_E\mid\mathcal B_{\mathcal N}]
  \mathbb E[\mathbf1_F\mid\mathcal B_{\mathcal N}]\\
 \hfill
 (E\in\mathcal B_{X\mathcal N},\ F\in\mathcal B_{Y\mathcal N}),\\[2pt]
 \mu_{(\bar M,X,Y)}\bigl(x\ind^K_{M_z}y\bigr)=1.
 \end{cases}
\]
In the last condition, $z$ denotes the restriction to $\bar m$ of
a type in $S_{\bar m,x,y}(T)$, and the locus is taken over
$z\in\mathcal M$. The locus is Borel by
Lemma~\ref{lem:app-borel-kim}.

If $\bar M'$ is another tuple as in
Lemma~\ref{lem:app-model-code}, choose, for each $i$, coordinates
$M'_{j(i,k)}$ of the same sort as $M_i$ such that
\[
 \sum_{k<\omega}d(M_i,M'_{j(i,k)})<\infty.
\]
In the joint distribution of $(\bar M,\bar M')$, in classical
variables $(\bar m,\bar m')$, Borel--Cantelli gives
\[
 m_i=m'_{j(i,k)}\text{ for all sufficiently large }k
 \quad\text{almost surely}.
\]
Applying the same argument with $\bar M,\bar M'$ interchanged
and intersecting the countably many conull sets gives, sort by sort,
\[
 \{m_i:i<\omega\}=\{m'_j:j<\omega\}
 \quad\text{almost surely}.
\]
The same approximations make each
$\llbracket\varphi(X,\bar M)\rrbracket$ a limit in measure of
events defined using $X,\bar M'$, and conversely. Thus
$\mathcal B_{X\mathcal N}$ is also independent of the
enumeration. It follows that $\ind^*$ is well defined and
invariant under automorphisms. Both defining conditions are
symmetric, so $\ind^*$ is symmetric.

\begin{lemma}\label{lem:app-random-coheir}
If $\tp^R(X/\mathcal NY)$ is finitely satisfiable in $\mathcal N$,
then $X\ind^*_{\mathcal N} Y$.
\end{lemma}

\begin{proof}
\begin{clm}
The measure $\mu_{(\bar M,X,Y)}$ gives measure one to the
locus with $z\in\mathcal M$ on which $\tp(x/M_z y)$ is
finitely satisfiable in $M_z$.
\end{clm}

\begin{clmproof}
Fix a formula $\varphi(x_0;y_0,\bar m_0)$ in finite subtuples of
$x,y,\bar m$, with corresponding random subtuples
$X_0,Y_0,\bar M_0$. Enumerate the tuples of coordinates of $\bar M$
of the sorts of $X_0$ by $(d_j)_{j<\omega}$. Put
\[
 F:=\llbracket\varphi(X_0;Y_0,\bar M_0)\rrbracket,\qquad
 U_n:=\bigcup_{j<n}\llbracket\varphi(d_j;Y_0,\bar M_0)\rrbracket,
 \qquad U:=\bigcup_nU_n.
\]
For tuples $C,D\in\mathcal N$ of these sorts,
\[
 \mathbb P\bigl(
  \llbracket\varphi(C;Y_0,\bar M_0)\rrbracket
  \mathbin{\triangle}
  \llbracket\varphi(D;Y_0,\bar M_0)\rrbracket
 \bigr)
 \leqslant\sum_{i<|X_0|}d(C_i,D_i).
\]
By density of $(d_j)_{j<\omega}$, for every such $C$,
\[
 \mathbb P\bigl(
 \llbracket\varphi(C;Y_0,\bar M_0)\rrbracket\setminus U
 \bigr)=0,
\]
and hence
\[
 \mathbb P\bigl(
 \llbracket\varphi(C;Y_0,\bar M_0)\rrbracket\setminus U_n
 \bigr)
 \leqslant\mathbb P(U\setminus U_n).
\]
For fixed $n$, the left-hand side is a continuous formula in
$C$ over $\mathcal NY$. Finite satisfiability transfers the bound
to $X_0$ with any positive error; letting that error tend to zero
gives
\[
 \mathbb P(F\setminus U_n)
 \leqslant\mathbb P(U\setminus U_n)\longrightarrow0,
 \qquad
 \mathbb P(F\setminus U)=0.
\]
There are countably many choices of
$\varphi,X_0,Y_0,\bar M_0$. Intersecting the corresponding
conull sets proves the claim.
\end{clmproof}

\begin{clm}
The algebras $\mathcal B_{X\mathcal N}$ and
$\mathcal B_{Y\mathcal N}$ are conditionally independent over
$\mathcal B_{\mathcal N}$.
\end{clm}

\begin{clmproof}
For $\CL$-formulas $\varphi,\psi$ in finite subtuples, put
\[
 f(X):=\mathbf1_{\llbracket\varphi(X,\bar M)\rrbracket},
 \qquad
 g(Y):=\mathbf1_{\llbracket\psi(Y,\bar M)\rrbracket},
 \qquad
 h:=\mathbb E[g(Y)\mid\mathcal B_{\mathcal N}].
\]
Choose simple functions
\[
 h_n=\sum_{i<k_n}r_{n,i}\mathbf1_{E_{n,i}},
 \qquad
 r_{n,i}\in\mathbb Q\cap[0,1],\quad
 E_{n,i}\in\mathcal B_{\mathcal N},\quad
 \|h_n-h\|_1\longrightarrow0,
\]
where $(E_{n,i})_{i<k_n}$ is a partition for each $n$.
For every $E\in\mathcal B_{\mathcal N}$ and every tuple
$C\in\mathcal N$ of the appropriate sorts, $\mathbf1_Ef(C)$
is $\mathcal B_{\mathcal N}$-measurable, so
\[
 \mathbb E[\mathbf1_Ef(C)g(Y)]
 =\mathbb E[\mathbf1_Ef(C)h],
\]
and
\[
 \left|\mathbb E[\mathbf1_Ef(C)g(Y)]
       -\mathbb E[\mathbf1_Ef(C)h_n]\right|
 \leqslant\|h-h_n\|_1.
\]
For fixed $E,n$, the left-hand side is a continuous formula
over $\mathcal NY$. Finite satisfiability and passage to the
limit give
\[
 \mathbb E[\mathbf1_Ef(X)g(Y)]
 =\mathbb E[\mathbf1_Ef(X)h]
 \qquad(E\in\mathcal B_{\mathcal N}).
\]
Therefore
\[
 \mathbb E[f(X)g(Y)\mid\mathcal B_{\mathcal N}]
 =\mathbb E[f(X)\mid\mathcal B_{\mathcal N}]
  \mathbb E[g(Y)\mid\mathcal B_{\mathcal N}].
\]
The monotone class theorem extends this identity to
$\mathcal B_{X\mathcal N}$ and $\mathcal B_{Y\mathcal N}$.
\end{clmproof}

The first claim implies
$\mu_{(\bar M,X,Y)}\bigl(x\ind^K_{M_z}y\bigr)=1$.
Together with the second claim, this gives
$X\ind^*_{\mathcal N} Y$.
\end{proof}

\subsection{Measurable independent amalgamation}

\begin{proposition}\label{prop:app-random-amalgamation}
Suppose that
\[
 X\ind^*_{\mathcal N} Y,\qquad
 X'\ind^*_{\mathcal N} W,\qquad
 Y\ind^*_{\mathcal N} W,\qquad X\equiv_{\mathcal N}X'.
\]
Then there is $X''$ such that
\[
 X''\equiv_{\mathcal N Y}X,\qquad
 X''\equiv_{\mathcal N W}X',\qquad
 X''\ind^*_{\mathcal N} YW.
\]
\end{proposition}

\begin{proof}
Let $x,y,w$ be the tuples of classical variables corresponding to
$X,Y,W$, respectively. The tuple $X'$ also corresponds to $x$.
Put
\[
 \lambda_{XY}:=\mu_{(\bar M,X,Y)},\qquad
 \lambda_{XW}:=\mu_{(\bar M,X',W)},\qquad
 \lambda_{YW}:=\mu_{(\bar M,Y,W)}.
\]
Disintegrate $\mu_{(\bar M,X)}$, $\mu_{(\bar M,Y)}$ and
$\mu_{(\bar M,W)}$ over $\kappa$, obtaining probability kernels
$\alpha_z(dp)$, $\beta_z(dq)$ and $\gamma_z(dr)$ on
$S_{\bar m,x}(T)$, $S_{\bar m,y}(T)$ and $S_{\bar m,w}(T)$,
respectively. These kernels concentrate on types extending $z$
for $\kappa$-almost every $z$. Since $X\equiv_{\mathcal N}X'$,
the kernel $\alpha_z$ also disintegrates $\mu_{(\bar M,X')}$.

Define the restriction maps
\[
 \pi_{XY}(s)=(s|\bar m x,s|\bar m y),\qquad
 \pi_{XW}(u)=(u|\bar m x,u|\bar m w).
\]
Conditional independence gives
\[
 \begin{split}
  (\pi_{XY})_*\lambda_{XY}
    &=\int\alpha_z\otimes\beta_z\,d\kappa(z),\\
  (\pi_{XW})_*\lambda_{XW}
    &=\int\alpha_z\otimes\gamma_z\,d\kappa(z).
 \end{split}
\]
Disintegrate $\lambda_{XY}$ and $\lambda_{XW}$ along these maps,
and denote the resulting probability kernels by
$K_{XY}(p,q;ds)$ and $K_{XW}(p,r;du)$.
Also write
\[
 \lambda_{YW}=\int D_z\,d\kappa(z),\qquad
 t\in S_{\bar m,y,w}(T),
\]
where $D_z$ concentrates on types extending $z$ for
$\kappa$-almost every $z$. Put $t_y=t|\bar m y$ and
$t_w=t|\bar m w$. Uniqueness of disintegration gives
\[
 (t\mapsto t_y)_*D_z=\beta_z,\qquad
 (t\mapsto t_w)_*D_z=\gamma_z
 \quad\text{for $\kappa$-almost every $z$}.
\]
All these disintegrations take place on standard Borel spaces;
the kernels concentrate on the corresponding fibres almost
everywhere.

Choose $z$ according to $\kappa$. Conditional on $z$, choose
the full type $t$ according to $D_z$ and, independently of $t$,
choose $p$ according to $\alpha_z$. Then choose $s$ and $u$
according to $K_{XY}(p,t_y;ds)$ and $K_{XW}(p,t_w;du)$,
independently conditional on $(z,t,p)$.
Retaining the full type $t$ preserves the prescribed conditional
distribution of $(\bar M,Y,W)$.
Thus let $\rho$ be the probability measure on tuples $(z,t,p,s,u)$
defined by
\begin{equation}\label{eq:app-configuration-law}
 \rho(dz,dt,dp,ds,du)
 =\kappa(dz)D_z(dt)\alpha_z(dp)
       K_{XY}(p,t_y;ds)K_{XW}(p,t_w;du).
\end{equation}

\begin{clm}
The $s$-, $u$- and $t$-restrictions of $\rho$ are
$\lambda_{XY}$, $\lambda_{XW}$ and $\lambda_{YW}$, respectively.
Moreover, $\rho$ gives measure one to the set of tuples satisfying
\[
 \begin{gathered}
  z\in\mathcal M,\qquad t|\bar m=p|\bar m=z,\\
  s|\bar m x=u|\bar m x=p,\qquad
  s|\bar m y=t_y,\qquad u|\bar m w=t_w,\\
  x\ind^K_{M_z} y\text{ in }s,\qquad
  x\ind^K_{M_z} w\text{ in }u,\qquad
  y\ind^K_{M_z} w\text{ in }t.
 \end{gathered}
\]
This set is Borel.
\end{clm}

\begin{clmproof}
For every bounded Borel function $f$ on $S_{\bar m,x,y}(T)$,
\[
 \begin{split}
  \int f(s)\,d\rho
   &=\int\kappa(dz)\int\alpha_z(dp)\int\beta_z(dq)
         \int f(s)K_{XY}(p,q;ds)\\
   &=\int f\,d\lambda_{XY}.
 \end{split}
\]
The calculation for $u$ is identical, and the assertion for $t$
follows directly from \eqref{eq:app-configuration-law}.
Concentration of the kernels on their fibres gives the displayed
restriction identities. The three Kim-independence conditions
hold almost surely by the hypotheses and the restriction
identities. Equality of restrictions is closed, $\mathcal M$ is
Borel, and Lemma~\ref{lem:app-borel-kim} gives Borelness of the
Kim-independence conditions.
\end{clmproof}

Denote the Borel conull set in the claim by $H$.

\begin{clm}
There is a Borel probability measure $\lambda$ on
$S_{\bar m,x,y,w}(T)$ such that, writing $z=v|\bar m$,
\[
 \begin{gathered}
  (v\mapsto v|\bar m xy)_*\lambda=\lambda_{XY},\qquad
  (v\mapsto v|\bar m xw)_*\lambda=\lambda_{XW},\\
  (v\mapsto v|\bar m yw)_*\lambda=\lambda_{YW},\qquad
  \lambda\bigl(x\ind^K_{M_z} yw\bigr)=1,
 \end{gathered}
\]
and the joint law of $(v|\bar m,v|\bar m x,v|\bar m yw)$ is
\[
                    \kappa(dz)\alpha_z(dp)D_z(dt).
\]
\end{clm}

\begin{clmproof}
Fix $(z,t,p,s,u)\in H$. Realize $t$ by $(\bar m,b,c)$,
realize $s$ over $\bar m b$ by $a_0$, and realize $u$ over
$\bar m c$ by $a_1$. Then
\[
 a_0\ind^K_{M_z} b,\qquad
 a_1\ind^K_{M_z} c,\qquad
 b\ind^K_{M_z} c,\qquad
 a_0\equiv_{M_z}a_1.
\]
By the independence theorem
\cite[Theorem~6.5]{KaplanRamsey2020}, there is
$v\in S_{\bar m,x,y,w}(T)$ with
\[
 v|\bar m xy=s,\qquad v|\bar m xw=u,\qquad
 v|\bar m yw=t,\qquad x\ind^K_{M_z} yw\text{ in }v.
\]
The relation on $H\times S_{\bar m,x,y,w}(T)$ defined by these
conditions is Borel, by Lemma~\ref{lem:app-borel-kim}, and all
its sections over $H$ are nonempty. Jankov--von Neumann
uniformization \cite[Theorem~29.9]{Kechris1995} gives a
universally measurable selection $F$ on $H$.
Fix $v_0\in S_{\bar m,x,y,w}(T)$ and put $F=v_0$ off $H$.
Since $H$ is Borel, this extension is universally measurable.
Put
\[
 \lambda(E):=\bar\rho\bigl(F^{-1}(E)\bigr)
 \qquad\bigl(E\subseteq S_{\bar m,x,y,w}(T)\text{ Borel}\bigr),
\]
where $\bar\rho$ is the completion of $\rho$.
Universal measurability makes $\lambda$ a well-defined Borel
probability measure. The restriction identities and
Kim-independence follow from the defining relation and the first
claim. Since
\[
 (F(z,t,p,s,u)|\bar m,
   F(z,t,p,s,u)|\bar m x,
   F(z,t,p,s,u)|\bar m yw)=(z,p,t)
 \quad\text{on }H,
\]
the last assertion follows from \eqref{eq:app-configuration-law}.
\end{clmproof}

Realize $\lambda$ in a sufficiently saturated randomization.
Its $(\bar m,y,w)$-restriction is $\mu_{(\bar M,Y,W)}$.
Quantifier elimination and homogeneity therefore allow us to use
the original tuple $(\bar M,Y,W)$ in this realization. Let $X''$
be the tuple corresponding to $x$. The other two restrictions give
\[
 X''\equiv_{\bar M Y}X,\qquad
 X''\equiv_{\bar M W}X'.
\]
These type equalities extend to $\mathcal N Y$ and
$\mathcal N W$, respectively: the coordinates of $\bar M$ are
dense in each random-element sort of $\mathcal N$, and their
formula events generate $\mathcal B_{\mathcal N}$.

Finally, let
\[
 E=\llbracket\varphi(\bar M,X'')\rrbracket,\qquad
 F=\llbracket\psi(\bar M,Y,W)\rrbracket
 \qquad(\varphi,\psi\in\CL).
\]
By Lemma~\ref{lem:app-model-code},
$\mathcal B_{\mathcal N}$ is the completed algebra generated
by the formula events of $\bar M$. Thus conditional
expectations over $\mathcal B_{\mathcal N}$ are represented
by measurable functions of $z$.
The joint restriction $\kappa(dz)\alpha_z(dp)D_z(dt)$ gives
\[
 \begin{aligned}
 \mathbb E[\mathbf1_E\mathbf1_F
                  \mid\mathcal B_{\mathcal N}](z)
 &=\alpha_z([\varphi])D_z([\psi])\\
 &=\mathbb E[\mathbf1_E\mid\mathcal B_{\mathcal N}](z)\,
   \mathbb E[\mathbf1_F\mid\mathcal B_{\mathcal N}](z)
 \end{aligned}
\]
for $\kappa$-almost every $z$.
A monotone-class argument extends this equality to
$E\in\mathcal B_{X''\mathcal N}$ and
$F\in\mathcal B_{YW\mathcal N}$; passing to the completed
algebras does not change the conditional expectations.
Together with $\lambda(x\ind^K_{M_z} yw)=1$, this proves
$X''\ind^*_{\mathcal N} YW$.
\end{proof}

\subsection{Proof of the theorem}

\begin{proof}[Proof of Theorem~\ref{thm:randomization-preserves-nsop1}]
Assume first that $T$ is countable and $\NSOP$, and that
some sort $\mathsf S$ has at least two elements.
The relation $\ind^*$ is invariant and symmetric.
Lemma~\ref{lem:app-random-coheir} gives condition (i) of
Lemma~\ref{lem:app-continuous-criterion}, and
Proposition~\ref{prop:app-random-amalgamation} gives condition
(iii) for random-element tuples.

We extend the relation to tuples with event-sort coordinates.
By fullness and completeness, choose
$u,v\in\mathsf S^R(\mathcal N)$ such that
\[
 \mathbb P\llbracket u\ne v\rrbracket=1,
\]
and include them among the coordinates of $\bar M$.
For an event $E$, let $w_E\in\mathsf S^R$ be determined by
\[
 \llbracket w_E=u\rrbracket=E,
 \qquad \llbracket w_E=v\rrbracket=E^c.
\]
The maps
\[
 E\longmapsto w_E,
 \qquad
 w\longmapsto\llbracket w=u\rrbracket
\]
are inverse definable maps over $u,v$ between the event sort and
\[
 \{w\in\mathsf S^R:
       \mathbb P\llbracket w=u\vee w=v\rrbracket=1\}.
\]
They preserve types and approximate finite satisfiability over
$\mathcal N$. Replace each event coordinate $E$ by $w_E$ and
apply $\ind^*$ to the resulting random-element tuples.
Changing $u,v$ gives codes which are interdefinable over $M_z$
almost surely and which generate the same event algebras over
$\mathcal B_{\mathcal N}$. Hence the extended relation is
well-defined and invariant, and it satisfies (i)--(iii) of
Lemma~\ref{lem:app-continuous-criterion}. That lemma gives that
$T^R$ is $\NSOP$.
If every sort of $T$ is a singleton, all random-element sorts
are trivial. With the atomless event sort retained in the
randomization, $T^R$ is therefore stable.

Now let $\CL$ be arbitrary. A continuous $\SOP$ witness in
$T^R$ uses only countably many symbols of $\CL$. Choose a
countable sublanguage $\CL_0$ containing them, together with
a non-singleton sort of $T$ if one exists, and put
$T_0=T|\CL_0$. Since $\SOP$ in a reduct implies $\SOP$ in the
original theory, $T_0$ is $\NSOP$. Randomization commutes with
this reduct, so the same formula and parameter tree would witness
$\SOP$ in $T_0^R$, contradicting the countable case.

Conversely, suppose that a classical formula $\theta(x;y)$ and
a tree $(b_\eta)_{\eta\in2^{<\omega}}$ witness $\SOP$ in $T$.
Let $B_\eta$ be the constant random tuple corresponding to $b_\eta$,
and let $X,Y$ be the random-element variables corresponding to
$x,y$. Put
\[
 Q(X;Y)=1-\mathbb P\llbracket\theta(X;Y)\rrbracket.
\]
For every $\sigma\in2^\omega$, a classical realization of
$\{\theta(x;b_{\sigma|n}):n<\omega\}$ gives a constant random
tuple realizing
\[
 \{Q(X;B_{\sigma|n})=0:n<\omega\}.
\]
Whenever $\eta\frown0\unlhd\xi$, the events
$\llbracket\theta(X;B_{\eta\frown1})\rrbracket$ and
$\llbracket\theta(X;B_\xi)\rrbracket$ are disjoint.
Consequently,
\[
 Q(X;B_{\eta\frown1})+Q(X;B_\xi)\geq1,
 \qquad
 \inf_X\max\{Q(X;B_{\eta\frown1}),Q(X;B_\xi)\}\geq\tfrac12.
\]
Thus $Q$ and the tree $(B_\eta)_{\eta\in2^{<\omega}}$ witness
$\SOP$ in $T^R$.
\end{proof}

\printbibliography
\end{document}